\documentclass[a4paper,12pt]{article}

\title{Improved bounds for the extremal number of subdivisions}
\date{}

\usepackage{pdfpages}
\usepackage{graphicx}
\usepackage{mathtools}
\usepackage{verbatim}
\usepackage{amssymb}
\usepackage{bm}
\usepackage{enumitem}
\usepackage{cases}
\usepackage{latexsym}
\usepackage{amsthm}
\usepackage{amsfonts}
\usepackage{bbm}
\usepackage{etoolbox}
\usepackage[top=25mm, bottom=25mm, left=28mm, right=28mm]{geometry}
\usepackage{cite}
\usepackage{changepage}

\newtheorem{theorem}{Theorem}
\AfterEndEnvironment{theorem}{\noindent\ignorespaces}
\newtheorem{corollary}[theorem]{Corollary}
\AfterEndEnvironment{corollary}{\noindent\ignorespaces}
\newtheorem{lemma}[theorem]{Lemma}
\AfterEndEnvironment{lemma}{\noindent\ignorespaces}

\AfterEndEnvironment{proposition}{\noindent\ignorespaces}

\AfterEndEnvironment{question}{\noindent\ignorespaces}

\theoremstyle{definition}
\newtheorem{definition}[theorem]{Definition}
\AfterEndEnvironment{definition}{\noindent\ignorespaces}

\theoremstyle{remark}

\AfterEndEnvironment{example}{\noindent\ignorespaces}

\AfterEndEnvironment{remark}{\noindent\ignorespaces}

\AfterEndEnvironment{notation}{\noindent\ignorespaces}

\newenvironment{proof2}{\proof[\textnormal{\textbf{Proof.}}]}{\qed}

\usepackage{chngcntr}
\usepackage{apptools}
\AtAppendix{\counterwithin{theorem}{section}}

\def\a{\alpha}

\def\d{\delta}

\def\s{\subset}
\def\ex{\text{ex}}

\begin{document}

\author{Oliver Janzer\thanks{Department of Pure Mathematics and Mathematical Statistics, University of Cambridge. E-mail: oj224@cam.ac.uk}}

\maketitle

\begin{abstract}

Let $H_t$ be the subdivision of $K_t$. Very recently, Conlon and Lee have proved that for any integer $t\geq 3$, there exists a constant $C$ such that $\ex(n,H_t)\leq Cn^{3/2-1/6^t}$. In this paper, we prove that there exists a constant $C'$ such that $\ex(n,H_t)\leq C'n^{3/2-\frac{1}{4t-6}}$.

\end{abstract}

\section{Introduction}

For a graph $H$, the extremal function $\ex(n,H)$ is defined to be the maximal number of edges in an $H$-free graph on $n$ vertices. This function is well understood for graphs $H$ with chromatic number at least three by the Erd\H os-Stone-Simonovits theorem. \cite{erdosstone,erdossimonovits} However, for bipartite graphs $H$, much less is known. For a survey on the subject, see \cite{furedisimonovits}. One of the few general results, proved by F\"uredi \cite{furedi}, and reproved by Alon, Krivelevich and Sudakov \cite{aks} is the following.

\begin{theorem}[F\"uredi, Alon-Krivelevich-Sudakov] \label{faks}
	Let $H$ be a bipartite graph such that in one of the parts all the degrees are at most $r$. Then there exists a constant $C$ such that $\ex(n,H)\leq Cn^{2-1/r}$.
\end{theorem}

Conlon and Lee \cite{conlonlee} have conjectured that the only case when this is tight up to the implied constant is when $H$ contains a $K_{r,r}$ (it is conjectured \cite{kst} that $\ex(n,K_{r,r})=\Omega(n^{2-1/r})$), and that for other graphs $H$ there exists some $\d>0$ such that $\ex(n,H)=O(n^{2-1/r-\d})$.

The subdivision of a graph $L$ is the bipartite graph with parts $V(L)$ and $E(L)$ (the vertex set and the edge set of the graph $L$, respectively) where $v\in V(L)$ is joined to $e\in E(L)$ if $v$ is an endpoint of $e$. It is easy to see that any $C_4$-free bipartite graph in which every vertex in one part has degree at most two is a subgraph of $H_t$ for some positive integer $t$, where $H_t$ is the subdivision of $K_t$. Conlon and Lee have verified their conjecture in the $r=2$ case by proving the following result.

\begin{theorem}[Conlon and Lee {\cite[Theorem 5.1]{conlonlee}}]
	For any integer $t\geq 3$, there exists a constant $C_t$ such that $\ex(n,H_t)\leq C_tn^{3/2-1/6^t}$.
\end{theorem}

They have observed the lower bound $\ex(n,H_t)\geq c_t n^{3/2-\frac{t-3/2}{t^2-t-1}}$ coming from the probabilistic deletion method, and have asked for an upper bound of the form $\ex(n,H)\leq C_tn^{3/2-\d_t}$, where $1/\d_t$ is bounded by a polynomial in $t$. We can prove such a bound even for a linear $\d_t$.

\begin{theorem} \label{subcomplete}
	For any integer $t\geq 3$, there exists a constant $C_t$ such that $\ex(n,H_t)\leq C_tn^{1+\frac{t-2}{2t-3}}=C_tn^{3/2-\frac{1}{4t-6}}$.
\end{theorem}

It would be very interesting to know whether or not this bound is tight up to the implied constant. It certainly is tight for $t=3$ as $\ex(n,C_6)=\Theta(n^{4/3})$.

We can in fact prove a slightly stronger result. For integers $s\geq 1$ and $t\geq 3$, let $L_{s,t}$ be the graph which is a $K_{s+t-1}$ with the edges of a $K_s$ removed. That is, the vertex set of $L_{s,t}$ is $S\cup T$ where $S\cap T=\emptyset$, $|S|=s$ and $|T|=t-1$, and $xy$ is an edge if and only if $x\in T$ or $y\in T$. Let $L'_{s,t}$ be the subdivision of $L_{s,t}$.

\begin{theorem} \label{mainresult}
	For any two integers $s\geq 1$ and $t\geq 3$, there exists a constant $C_{s,t}$ such that $\ex(n,L'_{s,t})\leq C_{s,t}n^{3/2-\frac{1}{4t-6}}$.
\end{theorem}

This result certainly implies Theorem \ref{subcomplete} as $L_{1,t}=K_t$. Moreover, we can apply Theorem \ref{mainresult} to obtain good bounds on the extremal number of the subdivision of the complete bipartite graph $K_{a,b}$ as well. Let us write $H_{a,b}$ for the subdivision of $K_{a,b}$. Conlon and Lee \cite[Theorem 4.2]{conlonlee} have proved that for any $2\leq a\leq b$ there exists a constant $C$ such that $\ex(n,H_{a,b})\leq Cn^{3/2-\frac{1}{12b}}$. They have also observed the lower bound $\ex(n,H_{a,b})=\Omega_{a,b}(n^{3/2-\frac{a+b-3/2}{2ab-1}})$ (which follows from the probabilistic deletion method). Hence their upper bound is reasonably close to best possible when $a=b$, but is weak when $b$ is much larger then $a$.

Since $K_{a,b}$ is a subgraph of $L_{b,a+1}$, Theorem \ref{mainresult} implies the following result, by taking $s=b$ and $t=a+1$.

\begin{corollary}
	For any two integers $2\leq a\leq b$, there exists a constant $C_{a,b}$, such $\ex(n,H_{a,b})\leq C_{a,b}n^{3/2-\frac{1}{4a-2}}$.
\end{corollary}

\section{Proof of Theorem \ref{mainresult}}

We shall use the following lemma of Conlon and Lee \cite[Lemma 2.3]{conlonlee}, which is a slight modification of a result of Erd\H os and Simonovits \cite{erdsimred}. Let us say that a graph $G$ is $\emph{K-almost-regular}$ if $\max_{v\in V(G)}\deg(v)\leq K\min_{v\in V(G)}\deg(v)$. Moreover, following Conlon and Lee, we say that a bipartite graph $G$ with a bipartition $A\cup B$ is \emph{balanced} if $\frac{1}{2}|B|\leq |A|\leq 2|B|$.

\begin{lemma}
	For any positive constant $\alpha<1$, there exists $n_0$ such that if $n\geq n_0$, $C\geq 1$ and $G$ is an $n$-vertex graph with at least $Cn^{1+\alpha}$ edges, then $G$ has a $K$-almost-regular balanced bipartite subgraph $G'$ with $m$ vertices such that $m\leq n^{\frac{\alpha(1-\alpha)}{2(1+\alpha)}}$, $|E(G')|\geq\frac{C}{10}m^{1+\a}$ and $K=60\cdot 2^{1+1/\a^2}$.
\end{lemma}

This reduces Theorem \ref{mainresult} to the following.

\begin{theorem} \label{mainreduced}
	For every $K\geq 1$, and positive integers $s\geq 1,t\geq 2$, there exists a constant $c=c(s,t,K)$ with the following property. Let $n$ be sufficiently large and let $G$ be a balanced bipartite graph with bipartition $A\cup B$, $|B|=n$ such that the degree of every vertex of $G$ is between $\d$ and $K\d$, for some $\d\geq cn^{\frac{t-2}{2t-3}}$. Then $G$ contains a copy of $L'_{s,t}$.
\end{theorem}

Given a bipartite graph $G$ with bipartition $A\cup B$, the \emph{neighbourhood graph} is the weighted graph $W_G$ on vertex set $A$ where the weight of the pair $uv$ is $d_G(u,v)=|N_G(u)\cap N_G(v)|$. Here and below $N_G(v)$ denotes the neighbourhood of the vertex $v$ in the graph $G$. For a subset $U\s A$, we write $W(U)$ for the total weight in $U$, ie. $W(U)=\sum_{uv\in {U \choose 2}} d_G(u,v)$.

We shall use the following simple lemma of Conlon and Lee \cite[Lemma 2.4]{conlonlee}.

\begin{lemma} \label{locallydense}
	Let $G$ be a bipartite graph with bipartition $A\cup B$, $|B|=n$, and minimum degree at least $\d$ on the vertices in $A$. Then for any subset $U\subset A$ with $\d|U|\geq 2n$,
	
	\begin{equation*}
		\sum_{uv\in {U \choose 2}} d_G(u,v)\geq \frac{\d^2}{2n}{|U| \choose 2}
	\end{equation*}	
\end{lemma}

\noindent In other words, the conclusion of Lemma \ref{locallydense} is that $W(U)\geq \frac{\d^2}{2n}{|U| \choose 2}$.

In the next definition, and in the rest of this paper, for a weighted graph $W$ on vertex set $A$, if $u,v\in A$, then $W(u,v)$ stands for the weight of $uv$. Moreover, we shall tacitly assume throughout the paper that $s\geq 1$ and $t\geq 3$ are fixed integers.

\begin{definition}
	Let $W$ be a weighted graph on vertex set $A$ and let $u,v\in A$ be distinct. We say that $uv$ is a \emph{light edge} if $1\leq W(u,v)<{s+t-1\choose 2}$ and that it is a \emph{heavy edge} if $W(u,v)\geq {s+t-1\choose 2}$.
\end{definition}

Note that if there is a $K_{s+t-1}$ in $W_G$ formed by heavy edges, then clearly there is an $L_{s,t}$ in $W_G$ formed by heavy edges, therefore there is an $L'_{s,t}$ in $G$.

The next lemma is one of our key observations.

\begin{lemma} \label{manylight}
	Let $G$ be an $L'_{s,t}$-free bipartite graph with bipartition $A\cup B$, $|B|=n$ and suppose that $W(A)\geq 8(s+t)^2n$. Then the number of light edges in $W_G$ is at least $\frac{W(A)}{4(s+t)^3}$.
\end{lemma}

\begin{proof2}
	Let $B=\{b_1,\dots,b_n\}$. Let $k_i=|N_G(b_i)|$ and suppose that $k_i\geq 2(s+t-2)$ for some $i$. As $G$ is $L'_{s,t}$-free, there is no $K_{s+t-1}$ in $W\lbrack N_G(b_i)\rbrack$ formed by heavy edges. Thus, by Tur\'an's theorem, the number of light edges in $N_G(b_i)$ is at least $$(s+t-2){\frac{k_i}{s+t-2} \choose 2}=\frac{1}{2}k_i\big(\frac{k_i}{s+t-2}-1\big)\geq \frac{k_i^2}{4(s+t-2)}.$$ But
	
	\begin{equation*}
		\sum_{i: k_i<2(s+t-2)} {k_i \choose 2}<4(s+t)^2n\leq \frac{W(A)}{2},
	\end{equation*} so
	
	\begin{equation*}
		\sum_{i: k_i\geq 2(s+t-2)} {k_i \choose 2}\geq \frac{W(A)}{2}.
	\end{equation*}
	Since every light edge is present in at most ${s+t-1 \choose 2}$ of the sets $N_G(b_i)$, it follows that the total number of light edges is at least
	
	$$ \frac{1}{{s+t-1 \choose 2}}\sum_{i: k_i\geq 2(s+t-2)} \frac{k_i^2}{4(s+t-2)}\geq \frac{W(A)}{4(s+t)^3}.$$
\end{proof2}

\begin{corollary} \label{lightcorollary}
	Let $G$ be an $L'_{s,t}$-free bipartite graph with bipartition $A\cup B$, $|B|=n$, and minimum degree at least $\d$ on the vertices in $A$. Then for any subset $U\subset A$ with $|U|\geq \frac{8(s+t)n}{\d}$ and $|U|\geq 2$, the number of light edges in $W_G\lbrack U\rbrack$ is at least $\frac{\d^2}{8(s+t)^3n}{|U| \choose 2}$.	
\end{corollary}

\begin{proof2}
	By Lemma \ref{locallydense}, we have $W(U)\geq \frac{\d^2}{2n}{|U| \choose 2}\geq \frac{\d^2}{8n}|U|^2\geq 8(s+t)^2n$. Now the result follows by applying Lemma \ref{manylight} to the graph $G\lbrack U\cup B\rbrack$. 
\end{proof2}

\medskip

We are now in a position to complete the proof of Theorem \ref{mainreduced}.

\begin{proof}[\textnormal{\textbf{Proof of Theorem \ref{mainreduced}}}]
	Let $c$ be specified later and suppose that $n$ is sufficiently large. Assume, for contradiction, that $G$ is $L'_{s,t}$-free. We shall find distinct $u_1,\dots,u_{t-1}\in A$ with the following properties.
	
	\begin{enumerate}[label=(\roman*)]
		\item Each $u_iu_j$ is a light edge in $W_G$
		\item If $i,j,k$ are distinct, then $N_G(u_i)\cap N_G(u_j)\cap N_G(u_k)=\emptyset$
		\item For each $1\leq i\leq t-1$, the number of $v\in A$ with the property that for every $j\leq i$, $u_jv$ is a light edge is at least $(\frac{\d^2}{32(s+t)^3n})^{i}\cdot|A|$
	\end{enumerate}

	
	As $n$ is sufficiently large, we have $|A|\geq n/2\geq \frac{8(s+t)n}{\d}$, therefore by Corollary \ref{lightcorollary} there are at least $\frac{\d^2}{8(s+t)^3n}{|A| \choose 2}$ light edges in $A$, so we may choose $u_1\in A$ such that the number of light edges $u_1v$ is at least $\frac{\d^2}{8(s+t)^3n}(|A|-1)\geq \frac{\d^2}{32(s+t)^3n}|A|$.
	
	Now suppose that $2\leq i\leq t-1$, and that  $u_1,\dots,u_{i-1}$ have been constructed satisfying (i),(ii) and (iii). Let $U_0$ be the set of vertices $v\in A$ with the property that $u_jv$ is a light edge for every $j\leq i-1$. By (iii), we have $|U_0|\geq (\frac{\d^2}{32(s+t)^3n})^{i-1}|A|$. Now let $U$ consist of those $v\in U_0$ for which $N_G(u_j)\cap N_G(u_k)\cap N_G(v)=\emptyset$ holds for all $1\leq j<k\leq i-1$. Since $u_ju_k$ is a light edge for any $1\leq j<k\leq i-1$, we have that $d_G(u_j,u_k)<{s+t-1 \choose 2}$. But the degree of every $b\in B$ is at most $K\d$, therefore the number of $v\in A$ for which $N_G(u_j)\cap N_G(u_k)\cap N_G(v)\neq \emptyset$ is at most ${s+t-1 \choose 2}K\d$, so $|U_0\setminus U|\leq {i-1 \choose 2}{s+t-1 \choose 2}K\d$. But note that for sufficiently large $n$, we have $(\frac{\d^2}{32(s+t)^3n})^{i-1}|A|\geq 2{i-1 \choose 2}{s+t-1 \choose 2}K\d$ because $\d=o((\d^2/n)^{t-2}n)$ and $\d=o((\d^2/n)n)$. Thus, 
	$$|U|\geq \frac{1}{2}|U_0|\geq \frac{1}{2}\big(\frac{\d^2}{32(s+t)^3n}\big)^{i-1}|A|.$$
	But for sufficiently large $c=c(s,t,K)$, we have $\frac{1}{2}(\frac{\d^2}{32(s+t)^3n})^{i-1}|A|\geq \frac{8(s+t)n}{\d}$. Indeed, this is obvious when $\d^2\geq 32(s+t)^3n$, and otherwise, using $\d\geq cn^{\frac{t-2}{2t-3}}$, we have
	$$\frac{1}{2}\big(\frac{\d^2}{32(s+t)^3n}\big)^{i-1}|A|\geq \frac{1}{2}\big(\frac{\d^2}{32(s+t)^3n}\big)^{t-2}|A|\geq \frac{1}{4(32(s+t)^3)^{t-2}}\cdot\frac{\d^{2t-4}}{n^{t-3}}\geq \frac{8(s+t)n}{\d}$$
	Thus, by Corollary \ref{lightcorollary}, there exists some $u_i\in U$ with at least $\frac{\d^2}{8(s+t)^3n}(|U|-1)\geq (\frac{\d^2}{32(s+t)^3n})^{i}|A|$ light edges adjacent to it in $U$. This completes the recursive construction of the vertices $\{u_j\}_{1\leq j\leq t-1}$.
	
	\medskip
	
	By (iii) for $i=t-1$, there is a set $V\s A$ consisting of at least $(\frac{\d^2}{32(s+t)^3n})^{t-1}|A|$ vertices $v$ such that for every $j\leq t-1$, $u_jv$ is a light edge.
	We shall now prove that there exist distinct $v_1,\dots,v_s\in V$ such that $N_G(u_i)\cap N_G(u_j)\cap N_G(v_k)\neq \emptyset$ for all $i\neq j$, and $N_G(u_i)\cap N_G(v_j)\cap N_G(v_k)\neq \emptyset$ for all $j\neq k$. It is easy to see that this suffices since then there is a copy of $L'_{s,t}$ in $G$, which is a subdivision of the copy of $L_{s,t}$ in $W_G$ whose vertices are $v_1,\dots,v_s,u_1,\dots,u_{t-1}$.
	
	\medskip
	
	We shall now choose $v_1,\dots,v_s$ one by one. Since every $u_iu_j$ is a light edge, the number of those $v\in A$ with $N_G(u_i)\cap N_G(u_j)\cap N_G(v)\neq \emptyset$ for some $i\neq j$ is at most ${t-1 \choose 2}{s+t-1 \choose 2}K\d$. Moreover, given any choices for $v_1,\dots,v_{k-1}\in V$, as each $u_iv_j$ is a light edge, the number of those $v\in A$ with $N_G(u_i)\cap N_G(v_j)\cap N_G(v)\neq \emptyset$ for some $i,j$ is at most $(t-1)(k-1){s+t-1 \choose 2}K\d$. Therefore as long as $|V|>{t-1 \choose 2}{s+t-1 \choose 2}K\d+(t-1)(s-1){s+t-1 \choose 2}K\d$, suitable choices for $v_1,\dots,v_s$ can be made. Since $|V|\geq (\frac{\d^2}{32(s+t)^3n})^{t-1}|A|$, this last inequality holds for large enough $c=c(s,t,K)$.	
\end{proof}

\bibliography{Bipturanbibliography}{}

\begin{thebibliography}{1}

\bibitem{aks}
Noga Alon, Michael Krivelevich, and Benny Sudakov.
\newblock Tur{\'a}n numbers of bipartite graphs and related {Ramsey}-type
  questions.
\newblock {\em Combinatorics, Probability and Computing}, 12:477--494, 2003.

\bibitem{conlonlee}
David Conlon and Joonkyung Lee.
\newblock On the extremal number of subdivisions.
\newblock {\em arXiv preprint arXiv:1807.05008}, 2018.

\bibitem{erdossimonovits}
Paul Erd{\H{o}}s and Mikl{\'o}s Simonovits.
\newblock A limit theorem in graph theory.
\newblock In {\em Studia Sci. Math. Hung}, 1965.

\bibitem{erdsimred}
Paul Erd{\H{o}}s and Mikl{\'o}s Simonovits.
\newblock Some extremal problems in graph theory.
\newblock In {\em Combinatorial theory and its applications}, 1969.

\bibitem{erdosstone}
Paul Erd{\H{o}}s and Arthur~H Stone.
\newblock On the structure of linear graphs.
\newblock {\em Bull. Amer. Math. Soc}, 52(1087-1091):1, 1946.

\bibitem{furedi}
Zolt{\'a}n F{\"u}redi.
\newblock On a {Tur{\'a}n} type problem of {Erd{\H o}s}.
\newblock {\em Combinatorica}, 11(1):75--79, 1991.

\bibitem{furedisimonovits}
Zolt{\'a}n F{\"u}redi and Mikl{\'o}s Simonovits.
\newblock The history of degenerate (bipartite) extremal graph problems.
\newblock In {\em Erd{\H{o}}s Centennial}, pages 169--264. Springer, 2013.

\bibitem{kst}
Tam{\'a}s K{\H o}v{\'a}ri, Vera S{\'o}s, and P{\'a}l Tur{\'a}n.
\newblock On a problem of {K}. {Zarankiewicz}.
\newblock In {\em Colloquium Mathematicum}, volume~1, pages 50--57, 1954.

\end{thebibliography}
\bibliographystyle{plain}

\end{document}